\documentclass[12pt,leqno]{article}
\usepackage{amsmath, amsthm, amsfonts, amssymb, color}
\usepackage{mathrsfs}
\theoremstyle{definition}

\newcommand{\scr}[1]{\mathscr #1}
\definecolor{wco}{rgb}{0.5,0.2,0.3}

\numberwithin{equation}{section} \theoremstyle{remark}

\newcommand{\ua}{\uparrow}

\title{{\bf General Extinction Results for   Stochastic Partial Differential Equations and Applications}\footnote{Supported in
 part by NNSFC(11131003), SRFDP, the Fundamental Research Funds for the Central Universities, and the DFG through CRC 701 and IRTG 1132.} }
\author{
{\bf  Michael R\"ockner$^{b)}$, Feng-Yu Wang$^{a),c)}$}\footnote{Corresponding author.
wangfy@bnu.edu.cn; F.Y.Wang@swansea.ac.uk}\\
\footnotesize{$^{a)}$School of Math. Sci. \& Lab. Math. Com. Sys.,
Beijing Normal
University, Beijing 100875, China}\\
\footnotesize{$^{b)}$Department of Mathematics, Bielefeld
University, D-33501 Bielefeld, Germany}\\
\footnotesize{$^{c)}$Department of Mathematics, Swansea University,
Singleton Park, SA2 8PP, UK} }
\begin{document}
\def\R{\mathbb R}  \def\ff{\frac} \def\ss{\sqrt} \def\B{\mathbf
B} \def\BB{\scr B}
\def\N{\mathbb N} \def\kk{\kappa} \def\m{{\bf m}}
\def\dd{\delta} \def\DD{\Delta} \def\vv{\varepsilon} \def\rr{\rho}
\def\<{\langle} \def\>{\rangle} \def\GG{\Gamma} \def\gg{\gamma}
  \def\nn{\nabla} \def\pp{\partial} \def\tt{\tilde}
\def\d{\text{\rm{d}}} \def\bb{\beta} \def\aa{\alpha} \def\D{\scr D}
\def\EE{\scr E} \def\si{\sigma} \def\ess{\text{\rm{ess}}}
\def\beg{\begin} \def\beq{\begin{equation}}  \def\F{\scr F}
\def\Ric{\text{\rm{Ric}}} \def\Hess{\text{\rm{Hess}}}
\def\e{\text{\rm{e}}} \def\ua{\underline a} \def\OO{\Omega}  \def\oo{\omega}
 \def\tt{\tilde} \def\Ric{\text{\rm{Ric}}}
\def\cut{\text{\rm{cut}}} \def\P{\mathbb P} \def\ifn{I_n(f^{\bigotimes n})}
\def\C{\scr C}      \def\aaa{\mathbf{r}}     \def\r{r}
\def\gap{\text{\rm{gap}}} \def\prr{\pi_{{\bf m},\varrho}}  \def\r{\mathbf r}
\def\Z{\mathbb Z} \def\vrr{\varrho} \def\ll{\lambda}
\def\L{\scr L}\def\Tt{\tt} \def\TT{\tt}\def\II{\mathbb I}
\def\i{{\rm i}}\def\Sect{{\rm Sect}}\def\E{\mathbb E} \def\H{\mathbb H}
\def\M{\scr M}\def\Q{\mathbb Q} \def\texto{\text{o}} \def\i{{\rm i}}
\def\O{\scr O}

\maketitle
\begin{abstract} Let $L$ be a positive definite self-adjoint operator on the $L^2$-space associated to a $\si$-finite measure space. Let
$H$ be the dual space of the domain of $L^{1/2}$ w.r.t. $L^2(\mu)$.  By using an It\^o type inequality for the $H$-norm and an integrability condition for the hyperbound of the semigroup $P_t:=\e^{-Lt}$, general extinction results are derived for a class of continuous adapted processes on $H$.
Main applications include stochastic and deterministic fast diffusion equations with fractional Laplacians. Furthermore,
we prove exponential integrability of the extinction time for all space dimensions in the singular diffusion version of the well-known Zhang-model for self-organized criticality,
provided the noise is small enough. Thus we obtain that the system goes to the critical state in finite time
 in the deterministic and with probability one in finite time in the stochastic case.

\end{abstract} \noindent
 AMS subject Classification:\ 60J75, 47D07.   \\
\noindent
 Keywords:   Extinction, Stochastic differential equation, fast diffusion equation.

 \vskip 2cm

\section{Introduction}
The phenomenon of self-organized criticality (SOC) is widely studied in physics from different perspectives.
In \cite{BJ92} it was proposed to describe this phenomenon, e.g. in the case of the avalanche dynamics in the
Bak-Tang-Wiesenfeld-model (see \cite{BTW88}) and in the Zhang-model (see \cite{Z89}) by a singular diffusion with the incoming
energy being realized by adding a noise term, for instance a linear multiplicative noise. The resulting stochastic partial diffusion equation
is then of the following type:
\beq\label{R1.1} \d X_t- \DD \Psi(X_t)\d t\ni B(X_t)\d W_t,\end{equation} where $\Psi:\R\to 2^\R$ is a multivalued map defined by
\beq\label{R1.1'} \Psi(s):= \beg{cases} \theta_1+\theta_2 s,\ &\text{if}\ s\in [0,\infty),\\
[0,\theta_1], \ &\text{if}\ s=0,\\
0, \ &\text{if}\ s<0,\end{cases}\end{equation}
where $\theta_1,\theta_2\ge 0$ are constants, with $\theta_2=0$ in the BTW-model whereas $\theta_2>0$ in the Zhang-model.
Hence in the BTW-model $\Psi$ is just the Heaviside function considered as a multivalued map, with jump at $s=0$, i.e.
we take the critical state to be equal to zero, which we may without loss of generality, by simply shifting the equation. In (\ref{R1.1}) we fix an open bounded domain $\O\subset \R^d$ and $\DD$ denotes the Dirichlet Laplacian on $\O$, $W_t$ is a cylindrical Wiener process
on $L^2(\O)$ with natural inner product $\<\cdot,\cdot\>_2$, but the solution process $X_t$ takes values in $H$ defined to be the completion of $L^2(\O)$ under the norm
$\|\cdot\|_H$ corresponding to the inner product $\<x,y\>_H:= \<x, (-\DD)^{-1}y\>_2,$ i.e. $H$ is just the dual of the classical
Sobolev space $H_0^{1,2}(\O)$ ( and  is usually denoted by $H^{-1}$). To explain the type of noise in (\ref{R1.1}), let $\{e_k\}_{k\ge 1}$ be
a normalized eigenbasis of $-\DD$ in $L^2(\O)$ with corresponding eigenvalues $\{\ll_k\}_{k\ge 1}$ numbered in increasing orders with multiplicities.
It is well-known that $\ll_k=\text{O}(k^{2/d})$ for large $k$. Then $B: H\to \scr L_2(L^2(\O),H)$ is defined by
\beq\label{R1.2} B(x)h= \sum_{k=1}^\infty \mu_k \<e_k, h\>_2x e_k,\ \ x\in H, h\in L^2(\O),\end{equation} for $\{\mu_k\}_{k\ge 1}\subset \R$
chosen in such a way that  $B(x)\in \scr L_2(L^2(\O), H), x\in H$, which is e.g. the case if there exists a constant $\vv\in (0,1)$ such that
(see (\ref{T}) below) $\sum_{k=1}^\infty \mu_k^2\ll_k^{\ff d 2\lor (1+\vv)}<\infty.$

The fundamental question about (\ref{R1.1}) is now whether the system will go to the critical state $(=0$ in our case) in finite
 time, i.e. letting
$$\tau_0:=\inf\{t\ge 0: \ |X(t)|_H=0\},$$
is $\tau_0<\infty$ for any initial condition $X_0=x\in H$, that is; do we have extinction in finite time?

Recently, it was proved in \cite{BDRb, BDPR11} (where \cite{BDPR11} is an improved version of \cite{BDRb}) that the answer
is yes if $d=1$, i.e. $\O\subset\R$, however, only with positive probability, that is, $\P(\tau_0<\infty)>0$, provided $X_0$ is not too far away
from $0$ (see \cite{BDRb, BDPR11} for details). We mention here that in \cite{BDRb} the $\mu_k, k\ge 1$, introduced above, were assumed to be zero starting from $k\ge N+1$, an  assumption that is dropped in \cite{BDPR11}. Furthermore, in both
\cite{BDRb, BDPR11} for simplicity $d\le 3$ was assumed. The question, however, whether
\beq\label{P1}\P(\tau_0<\infty)=1\end{equation} was
left open and it seemed to be out of reach  of the methods in \cite{BDRb, BDPR11}.

The purpose of this paper is twofold. First, we develop a general technique to prove extinction for solutions of stochastic (ordinary and partial)
differential equations, and second to apply the outcoming results to equations of type (\ref{R1.1}) (in fact for a whole class of $\Psi: \R\to 2^\R$)
to analyze the above question, more precisely the problem whether one of the following three properties hold:
\beg{enumerate} \item[(i)] $\P(\tau_0<\infty)>0$ for small $X_0,$
\item[(ii)] $\P(\tau_0<\infty)>0$ for all $X_0,$
\item[(iii)] $\P(\tau_0<\infty)=1$ or, moreover, $\tau_0$ has finite  (exponential or polynomial) moments.   \end{enumerate}

We already want to mention here that as a consequence we obtain that for all dimensions in the Zhang-model of SOC the extinction time $\tau_0$ is even exponentially integrable, in particular we have extinction in finite time with probability one, provided the noise is small enough and
$X_0=x\in L^4(\O), x\ge 0.$ This has been open even in the deterministic case. In order to also include stochastic fast diffusion equations and prove extinction for that case, strengthening and generalizing the results from \cite{BDRa} we study (\ref{R1.1}) for the following class of
$\Psi: \R\to 2^\R$: $\Psi$ is a maximal monotone graph with $0\in\Psi(0)$ such that for some constants $r\in [0,1), C>0, \theta_1>0, q\ge 1+r,$ and $\theta_2\ge 0$
\beq\label{R1.4} C(|s|^q+|s|)\ge sy \ge \theta_1|s|^{1+r}+\theta_2|s|^2,\ \ \ s\in\R, y\in \Psi(s).\end{equation} In this case we call
(\ref{R1.1}) fast diffusion-type equation. A special case of this is the stochastic fast diffusion equation, i.e. equation (\ref{R1.1}) with
\beq\label{R1.4'} \Psi(s)= s^{r-1}|s|,\ \ \ s\in\R,\end{equation} where $r\in (0,1)$. In this case we, however,
avoid the assumption $d\le 3$, made in \cite{BDRb} for simplicity, but prove extinction with positive probability  in all
dimensions with the usual dimension dependent restriction on $r$ known from the deterministic case.
 As in \cite{BDRa, BDRb, BDPR11} the latter is, of course,
always included in our results choosing $\mu_k=0$ for all $k\ge 1$ in (\ref{R1.2}).

Another new feature of this paper is that we prove our results for a whole class of general operators $L$ on a measurable space replacing the Dirichlet Laplacian $(-\DD)$ on $\O\subset \R^d$. This class, in particular, includes fractional Laplacians $L=(-\DD)^\aa, \aa\in (0,1]$,
which have attracted more interest recently. Existence and uniqueness of solutions to (\ref{R1.1}), with $(-\DD)$ replaced by a fractional Laplacian, have first been proved in \cite{RRW} in the general, i.e. stochastic,   hence including the deterministic $(B=0)$, case. In the deterministic case
these results have been reproved in \cite{V10} under some restrictions on dimensions, however, by  completely different methods.

The organization of this paper is as follows. In Section 2 we state and prove our two general results on extinction of solutions for
stochastic equations. The first (see Theorem \ref{T2.1} below) gives quantitative conditions so that properties (i)-(iii) above
 hold respectively. It also confirms the intuition that one gets stronger extinction results if one has more coercivity
 in the system (i.e. $\rr_2>0$ in condition $(H1)$ below or correspondingly $\theta_2>0$ in (\ref{R1.4}) above, as e.g. in the Zhang-model). The reason is that
 more coercivity means a stronger drift towards zero, so that the part in (\ref{R1.4}) with $\theta_1$ in front becomes very big and pushes the
 process to zero. The second result (see Theorem \ref{T2.2} below) tells us that suitable large enough noise has the same effect.
 In Section 3 these results are applied to the fast diffusion-type equations above, but for the said general class of operators
  replacing $(-\DD)$. Subsection 3.1 is devoted to the case, where $L= (-\DD+V)^\aa, \aa\in (0,1],$ and $V$ a nonnegative measurable
  function. In Subsection 3.2 we prove exponential integrability of the extinction time $\tau_0$ in the strongly dissipative case even for uncoloured
  noise (i.e. $\mu_k=1$ for all $k\ge 1$ in (\ref{R1.2})). In Subsection 3.3 we give examples of noises which lead to extinction with probability one.

  The above described SOC case, i.e. $\Psi$ is given by (\ref{R1.1'}), and the fast diffusion case (\ref{R1.4'}), both for the Laplacian and the fractional Laplacian, are considered as guiding examples, and are discussed in detail in Subsection 3.1.

\section{A general result}

Let $(E,\scr E, \mu)$ be a $\si$-finite measure space, and let $(L,\D(L))$ be a positive definite self-adjoint operator on $L^2(\mu)$ such that for its spectrum $\si(L)$ we have $\inf \si(L)>0$. Let $P_t=\e^{-Lt}$.
Let $H$ be the completion of $L^2(\mu)$  w.r.t. the norm $\|\cdot\|_H$ corresponding to the inner product $\<x,y\>_H:=\<x, L^{-1}y\>_2$ (note that since $\inf\si(L)>0$,
$L: \D(L)\to L^2(\mu)$ is bijective. So, this definition makes sense). For any $p,q\ge 1,$ let $\|\cdot\|_p$ and $\|\cdot\|_{p\to q}$  the norm on $L^p(\mu)$ and the operator norm from $L^p(\mu)$ to $L^q(\mu)$ respectively.

Let $\{X_t\}_{t\ge 0}$ be an $H$-valued continuous adapted process on the complete filtered probability space $(\OO,\{\F_t\}_{t\ge 0}, \P)$. Let
$$\tau_0:=\inf\{t\ge 0: X_t=0\}$$ be the extinction time of the process. To investigate the finiteness of $\tau_0$, we introduce the following conditions:
 \beg{enumerate} \item[$(H1)$] There exist $r\in [0,1)$ and constants $\rr_1>0, \rr_2,\rr_3\ge 0$ such that $\|X_t\|_H^2$ is a semimartingale,
$\|X_t\|_2^2$ and $\|X_t\|_{1+r}^{1+r}$ are locally integrable in $t$, and the It\^o differential of $\|X_t\|_H^2$ satisfies
$$\d\|X_t\|_H^2\le 2\big\{\rr_3\|X_t\|_H^2 -\rr_1\|X_t\|_{1+r}^{1+r}-\rr_2\|X_t\|_2^2\big\}\d t +\d M_t$$ for some local martingale $M_t$.
\item[$(H2)$] There exists $\theta\in (0,1]$ such that
$$\gg(\theta):= \int_0^\infty \e^{-\ll_1 (1-\theta)t}\|P_t\|_{1+r\to \ff{1+r}r}^\theta\d t<\infty,$$ where $\ll_1:=\inf\si(L)>0$ and $\si(L)$ is the spectrum of $L$. \end{enumerate}
Note that $(H1)$ implies that for any $X_0\in L^2(\OO\to H, \F_0, \P)$ and $T>0$, one has $(X_t)_{t\in [0,T]}\in L^{1+r}([0,T]\times \OO\to  L^{1+r}(\mu);\d t\times \P)$ and furthermore,
$(X_t)_{t\in [0,T]}\in L^{2}([0,T]\times \OO\to  L^{2}(\mu);\d t\times \P)$ if $\rr_2>0.$

\paragraph{Remark 2.1.}  $(H1)$ implies that $\e^{-2\rr_3 t} \|X_t\|_H^2, t\ge 0,$ is a non-negative local supermartingale. Hence it is equal to zero for all $t\ge\tau_0$ by a well-known result (see e.g. \cite[Chap IV, Lemma 3.19]{MR92}).

%We remark that for any $X_0=x\in H$,   $(H1)$ implies
%\beq\label{**} \E\bigg\{\ff{1_{\{\tau_0<\infty\}}\|X_{s+\tau_0}\|_H^2}{1+\|X_{s+\tau_0}\|_H^2} \bigg\}=0,\ \ \ s\ge 0,\end{equation}
%so that $X_t=0$ holds for $t\ge \tau_0.$ Thus, $\tau_0$ is the extinction time of $X_t$. To see this, we observe that by $(H1)$ the process
%$\xi_t:= \ff{\|X_t\|_H^2}{1+\|X_t\|_H^2}$ satisfies
%$$\d\xi_t \le c \xi_t\d t +\d\tt M_t$$ for some local martingale $\tt M_t$. Let $\tau_n\uparrow\infty$ be a sequence of stopping times such that for each $n\ge 1$, $(\tt M_{t\land \tau_n})_{t\ge 0}$ is a martingale. Then
%$$\E\big\{1_{\{\tau_0<\infty\}}\xi_{(s+\tau_0)\land\tau_n}\big\}\le c \int_0^s \E\big\{1_{\{\tau_0<\infty\}}\xi_{(u+\tau_0)\land\tau_n}\big\}\d u
%+\E\big\{1_{\{\tau_0<\infty\}}\big(\tt M_{(s+\tau_0)\land\tau_n}-\tt M_{\tau_0\land\tau_n}\big)\big\}.$$
%Noting that by the Doob stopping theorem
% $$\E\big\{1_{\{\tau_0<\infty\}}\big(\tt M_{(s+\tau_0)\land\tau_n}-\tt M_{\tau_0\land\tau_n}\big)\big\} =\E\big\{1_{\{\tau_0<\infty\}}\big(\E(\tt M_{(s+\tau_0)\land\tau_n}|\scr F_{\tau_0})-\tt M_{\tau_0\land\tau_n}\big)\big\}=0,$$  the Gronwall lemma yields that
 %$$\E\big\{1_{\{\tau_0<\infty\}}\xi_{(s+\tau_0)\land\tau_n}\big\}\le \e^{cs} \E\big\{1_{\{\tau_0<\infty\}}\xi_{\tau_0\land\tau_n}\big\}.$$
% Since $(\xi_t)_{t\ge 0}$ is bounded and $1_{\{\tau_0<\infty\}}\xi_{\tau_0}=0,$ by letting $n\to\infty$ and using the dominated convergence theorem we prove (\ref{**}).

\paragraph{Remark 2.2.} We would like to indicate that assumptions $(H1)$ and $(H2)$ are fulfilled for a large class of SPDEs. More precisely:
\beg{enumerate}\item[(1)] Let e.g. $L=-\DD+V$, where $\DD$ is the Dirichlet Laplacian on an open domain $\O\subset \R^d$ and $V\ge 0$ is a continuous function on $\O$. If either $\O$ is bounded or $\liminf_{|x|\to\infty}V(x)>0$, then $\ll_1>0$ and
$$\|P_t\|_{1+r\to \ff{1+r}r} \le c_1 \e^{-c_2 t} t^{-d(1-r)/[2(1+r)]},\ \ t>0$$ holds for some constants $c_1,c_2>0$ (see Section 3 for details). Therefore, $(H2)$ holds for $\theta \in (0, \ff{2(1+r)}{d (1-r)})\cap (0,1].$
\item[(2)] Let $X_t$ solve the following SPDE on $H$:
\beq\label{E1} \d X_t + L\Psi(X_t)\ni B(X_t)\d W_t,\end{equation} where $W_t$ is a cylindrical Brownian motion on $L^2(\mu)$; $\Psi: \R\to 2^\R$ is a maximal monotone graph such that
$$C(| s|^{1+r}+  |s|)\ge  sy\ge \theta_1 |s|^{1+r}+\theta_2 |s|^2,\ \ \  s\in\R, y\in\Psi(s)$$ for   some constants $C>0, \theta_1>0, \theta_2\ge 0$; and $B: H\to \scr L_2(L^2(\mu); H)$ (the space of Hilbert-Schmidt linear operators from $L^2(\mu)$ to $H$) such that $$\|B(x)\|_{\scr L_2(L^2(\mu);H)}^2\le 2\rr_3\|x\|_H^2,\ \ x\in H.$$ For what is meant by a solution to (\ref{E1}), we refer to Section 3 below.
Then, by  It\^o's  formula  $(H1)$ holds (see Section 3 below). \end{enumerate}

\beg{thm}\label{T2.1} Assume $(H1)$ and $(H2)$ and let $\E \|X_0\|_H^{\theta (1-r)}<\infty.$ \beg{enumerate}
\item[$(1)$] If either $\rr_2>0$ or $\rr_2=0$ but $\theta =1$, then
$$\P(\tau_0=\infty) \le 1-\E \e^{-\theta (1-r)\rr_3\tau_0} \le \ff{\rr_3 \gg(\theta)^{(1+r)/2}}{\rr_1^\theta \rr_2^{1-\theta} \ll_1^{(1-r)(1-\theta)/2}}\E \|X_0\|_H^{\theta (1-r)},$$ where $\rr_2^{1-\theta}:=1,$ if $\theta=1$ and $\rr_2=0.$ Consequently, $\P(\tau_0<\infty)>0$ holds for small $\E \|X_0\|_H^{\theta (1-r)}.$ \item[$(2)$] If $\rr_3<\ll_1 \rr_2$ then for any $\bb \in (0, \theta (1-r) (\rr_2\ll_1 -\rr_3)),$
$$\E\e^{\bb \tau_0}\le 1+ \ff{\bb \gg(\theta)^{(1+r)/2}\E\|X_0\|_H^{\theta(1-r)}}{\theta (1-r)\rr_1^\theta \aa_\bb^{1-\theta} \ll_1^{(1-r)(1-\theta)/2}}<\infty$$ holds for $\aa_\bb:= \rr_2 - \ff 1 {\ll_1} (\ff\bb {\theta(1-r)} +\rr_3)\in (0,\rr_2).$
\item[$(3)$] If $\rr_3= \ll_1\rr_2$ and $\theta=1$, then
$$\E\tau_0\le \ff{\gg(1)^{(1+r)/2}}{(1-r)\rr_1} \E\|X_0\|_H^{1-r}<\infty.$$\end{enumerate}\end{thm}

\paragraph{Remark 2.3.} If in $(H1)$ we have that $\rr_2=0$ and $M=0$, which is e.g. the case for deterministic stochastic fast
 diffusion equations (see equation (\ref{EE}) with $B\equiv 0$ in Section 3 below), it follows by (1) (or (3)) in Theorem \ref{T2.1} that we have extinction in finite time. Thus we recover the well-known results from the deterministic case.

\

Obviously, in cases (2) and (3) of Theorem \ref{T2.1} one has $\P(\tau_0<\infty)=1$ for all deterministic initial data $X_0=x\in H$,
while in case (1) $\P(\tau_0<\infty)>0$ holds for small enough $X_0=x\in H$. Our next result strengthens the assertion in (1):
 if the local martingale $M_t$ is strong enough but not too large,  then   $\P(\tau_0<\infty)>0$ holds  for all $X_0=x\in H$.

\beg{thm}\label{T2.2} Assume that $(H1)$ and $(H2)$ hold with either $\rr_2>0$, or $\rr_2=0$ but $\theta =1$. If there exist two functions  $g_2\ge g_1\in C([0,\infty))$ with $g_1(s)>0$ for $s>0$ such that
\beq\label{C}g_2(\|X_t\|_H^2)\d t\ge \d\<M\>_t\ge g_1(\|X_t\|_H^2)\d t.\end{equation} define
$$\xi(s):= \ff{\rr_1^\theta \rr_2^{1-\theta}\ll_1^{(1+r)(1-\theta)/2}}{\gg(\theta)^{(1+r)/2}}s^{1-\theta(1-r)/2} -\rr_3 s,\ \ s\ge 0,$$ where
$\rr_2^{1-\theta}:=1$ if $\rr_2=0$ and $\theta=1;$  and set
$$g= g_21_{\{\xi\ge 0\}} +g_11_{\{\xi<0\}}.$$ Then for any $X_0=x\in H$,
\beq\label{C'}\P(\tau_0=\infty)\le \ff{\int_0^{\|x\|_H^2}\exp\big[2\int_1^t \ff{\xi(u)}{g(u)}\d u\big]\d t }{\int_0^\infty \exp\big[2\int_1^t \ff{\xi(u)}{g(u)}\d u\big]\d t }<1.\end{equation}
 Consequently, if $g_1(s)\ge 2\rr_3s^2$ holds for large $s$, then $\P(\tau_0<\infty)=1.$ \end{thm}

 Theorem \ref{T2.2} tells us that if the quadratic variation of $M$ is big enough (but not too big, e.g. to avoid explosion in finite time),
 then we have extinction for all initial data in $H$. Intuitively, this means that in this case the process $X$ will come close
 enough to zero with positive probability,  so that by the strong Markov property Theorem \ref{T2.1} (1) applies.

 \

To prove Theorems \ref{T2.1} and \ref{T2.2}, we need the following two basic lemmas.

\beg{lem}\label{L2.3} $\|x\|_H\le \ss{\gg(\theta) }\, \|x\|_2^{1-\theta}\|x\|_{1+r}^\theta,\ \ \ x\in L^2(\mu)\cap L^{1+r}(\mu).$ \end{lem}
\beg{proof} Since by the symmetry of $P_t$ and since $\inf\si(L)=\ll_1,$ we have $$\mu(xP_t x)=\|P_{t/2}x\|_2^2\le \e^{-\ll_1 t}\|x\|_2^2,$$ and
by the H\"older inequality
$$\mu(xP_t x)\le \|x\|_{1+r} \|P_t x\|_{(1+r)/r}\le \|x\|_{1+r}^2 \|P_t\|_{1+r\to \ff{1+r}r},$$ it follows that
\beg{equation*}\beg{split}&\|x\|_H^2 =\int_0^\infty \mu(xP_tx)\d t\\
 &\le \|x\|_2^{2(1-\theta)}\|x\|_{1+r}^{2\theta} \int_0^\infty\e^{-\ll_1(1-\theta)t}\|P_t\|_{1+r\to\ff{1+r}r}^\theta\d t =\gg(\theta) \|x\|_2^{2(1-\theta)}\|x\|_{1+r}^{2\theta}.\end{split}\end{equation*} This completes the proof.\end{proof}

 \beg{lem}\label{L2.4} For any $\theta\in (0,1]$,
 $$b^{(1-\theta)/\theta} +\ff a b \ge a^{1-\theta},\ \ \ b>0, a\ge 0,$$ where $a^{1-\theta}:=1$ for $\theta=1.$ \end{lem}
\beg{proof} If $b\ge a^\theta$ then
$$b^{(1-\theta)/\theta} +\ff a b \ge  b^{(1-\theta)/\theta} \ge a^{1-\theta};$$ while if $b\le a^\theta$ then
$$b^{(1-\theta)/\theta} +\ff a b \ge  \ff a b \ge a^{1-\theta}.$$\end{proof}

\beg{proof}[Proof of Theorem \ref{T2.1}]
 By $(H1)$ and  It\^o's  formula, we have
\beq\label{I1} \beg{split}\d \|X_t\|_H^{\theta (1-r)} \le &-\theta (1-r) \|X_t\|_H^{\theta(1-r)-2} \big\{\rr_1\|X_t\|_{1+r}^{1+r}+\rr_2\|X_t\|_2^2 -\rr_3\|X_t\|_H^2\big\}\d t\\
 &+\d \tt M_t,\ \ t<\tau_0,\end{split}\end{equation} where $$\d \tt M_t = \theta(1-r) \|X_t\|_H^{\theta(1-r)-2}\d M_t,\ \ t<\tau_0.$$
Let $\aa\in (0,\rr_2]$  and $\aa=0$ if $\rr_2=0, \theta=1.$  By Lemma \ref{L2.3} and Lemma \ref{L2.4} for
$$b:= \ff{\|X_t\|_H^{1+r}}{\|X_t\|_2^{1+r}},\ \ a:= \ff{\ll_1^{(1-r)/2}\aa\gg(\theta)^{(1+r)/(2\theta)}\|X_t\|_H^{1-r}}{\rr_1},$$
and noting that $\|X_t\|_2^2\ge \ll_1 \|X_t\|_H^2$, we obtain for $t<\tau_0$
\beq\label{A0} \beg{split}&\rr_1 \|X_t\|_{1+r}^{1+r} +\rr_2 \|X_t\|_2^2 - \rr_3\|X_t\|_H^2\\
&\ge \ff{\rr_1\|X_t\|_H^{(1+r)/\theta}}{\gg(\theta)^{(1+r)/(2\theta)}  \|X_t\|_2^{(1-\theta)(1+r)/\theta}}+\aa\|X_t\|_2^2+\big\{(\rr_2-\aa)\ll_1 -\rr_3\big\}\|X_t\|_H^2\\
&\ge \ff{\rr_1\|X_t\|_H^{1+r}}{\gg(\theta)^{(1+r)/(2\theta)} }\bigg(b^{(1-\theta)/\theta}
 +\ff a b  \bigg)  +\big\{(\rr_2-\aa)\ll_1 -\rr_3\big\}\|X_t\|_H^2\\
&\ge \ff{\rr_1\|X_t\|_H^{1+r}}{\gg(\theta)^{(1+r)/(2\theta)} }\bigg(\ff{\ll_1^{(1-r)/2}\aa\gg(\theta)^{(1+r)/(2\theta)}\|X_t\|_H^{1-r}}{\rr_1}\bigg)^{1-\theta} +\big\{(\rr_2-\aa)\ll_1 -\rr_3\big\}\|X_t\|_H^2\\
&= \ff{\rr_1^\theta \aa^{1-\theta} \ll_1^{(1-r)(1-\theta)/2} \|X_t\|_H^{2-\theta(1-r)}}{\gg(\theta)^{(1+r)/2}}+\big\{(\rr_2-\aa)\ll_1 -\rr_3\big\}\|X_t\|_H^2.\end{split}\end{equation} Let
$$c_1:= \ff{\theta (1-r)\rr_1^\theta \aa^{1-\theta}\ll_1^{(1-r)(1-\theta)/2}}{\gg(\theta)^{(1+r)/2}}>0,\ \ c_2:= \theta(1-r) \big\{(\rr_2-\aa)\ll_1-\rr_3\big\}\in\R.$$ Combining (\ref{I1}) with (\ref{A0}) we obtain
$$\d\|X_t\|_H^{\theta (1-r)}\le -c_1\d t - c_2 \|X_t\|_H^{\theta (1-r)}\d t +\d\tt M_t,\ \ t<\tau_0.$$ Therefore,
$$\d\{\|X_t\|_H^{\theta (1-r)} \e^{c_2 t}\} \le -c_1\e^{c_2 t} \d t +\e^{c_2 t}\d\tt M_t, \  \ t<\tau_0.$$ This implies
$$ c_1 \E\int_0^{(\tau_0\land t-\vv)^+} \e^{c_2 s}\d s \le \E \|X_0\|_H^{\theta (1-r)},\ \ \ t,\vv>0.$$ Letting $\vv\to 0$ and $t\to\infty$ we arrive at
\beq\label{A1} \ff{\E \e^{c_2 \tau_0}-1}{c_2}= \E\int_0^{\tau_0} \e^{c_2 s}\d s \le \ff {\E\|X_0\|_H^{\theta(1-r)}}{c_1},\end{equation} where
$\ff{\E \e^{c_2 \tau_0}-1}{c_2}:=\E\tau_0$ if $c_2=0$. From this we are able to prove the desired assertions as follows.

(1) If either $\rr_2>0$ or $\rr_2=0$ but $\theta =1$, we take $\aa=\rr_2$. Then $c_1>0$ and $c_2=-\theta(1-r)\rr_3\le 0$. By (\ref{A1}) we obtain
$$\ff{\E\|X_0\|_H^{\theta(1-r)}}{c_1} \ge \ff{1-\E \e^{-\theta (1-r)\rr_3\tau_0}}{\theta (1-r)\rr_3},$$ which implies the second inequality by the definition of $c_1$. The first inequality is trivial if $\rr_3>0$. When $\rr_3=0$ we have $c_2=0$, so that (\ref{A1}) implies that $\E \tau_0<\infty$
and thus, the first inequality remains true.

(2) Let $\rr_3<\ll_1\rr_2$ and $\bb\in (0, \theta(1-r)(\ll_1\rr_2-\rr_3)).$ Take $\aa=\aa_\bb.$ We have $c_2=\bb>0$. So, (\ref{A1}) yields that
$$\E\e^{\bb \tau_0} \le 1 + \ff{\bb \E\|X_0\|_H^{\theta(1-r)}}{c_1}= 1+ \ff{\bb \gg(\theta)^{(1+r)/2}\E\|X_0\|_H^{\theta(1-r)}}{\theta (1-r)\rr_1^\theta \aa_\bb^{1-\theta} \ll_1^{(1-r)(1-\theta)/2}}.$$

(3) If $\rr_3=\ll_1\rr_2$ and $\theta=1$, we take $\aa=0$ so that $c_1>0$ and $c_2=0.$ Therefore, (\ref{A1}) implies that
$$\E \tau_0 \le \ff{\E\|X_0\|_H^{1-r}}{c_1} =\ff{\gg(1)^{(1+r)/2}}{(1-r)\rr_1} \E\|X_0\|_H^{1-r}.$$
\end{proof}

\beg{proof}[Proof of Theorem \ref{T2.2}] By (\ref{A0}) with $\aa=\rr_2$ and $(H1)$ we have for $t<\tau_0$
\beq\label{B1} \d\|X_t\|_H^2 \le -\xi(\|X_t\|_H^2) \d t +\d M_t.\end{equation} By the definition of $\xi$ we see that there exists a constant $r_0>0$ such that $\xi$ is strictly positive and   increasing  on $(0,r_0],$ and
\beq\label{BB2} \int_0^{r_0} \ff 1 {\xi(t)}\d t<\infty.\end{equation} For any constant $N> r_0\lor \|x\|_H^2,$ let
$$\tau_N:= \inf\{t\ge 0: \|X_t\|_H^2>N\}$$ and
$$f_N(s):= \int_0^s\d t \int_t^N\ff 2 {g_1(u)}\exp\bigg[-2\int_t^u\ff{\xi(v)}{g_1(v)}\,\d v\bigg]\d u,\ \ s\in (0,N].$$
Since $\xi$ is strictly positive and increasing on $(0,r_0]$, for $t\in (0,r_0]$ one has
\beg{equation*}\beg{split} &\int_t^{r_0}\ff 2 {g_1(u)}\exp\bigg[-2\int_t^u\ff{\xi(v)}{g_1(v)}\,\d v\bigg]\d u\\
&\le \ff 1 {\xi(t)}\int_t^{r_0}\ff {2 \xi(u)} {g_1(u)}\exp\bigg[-2\int_t^u\ff{\xi(v)}{g_1(v)}\,\d v\bigg]\d u\le \ff {1} {\xi(t)}.
\end{split}\end{equation*}   Combining this with (\ref{BB2}) we conclude that $f_N\in C^2((0,N]).$ Moreover, it is easy to check that $f_N''\le 0$ and $$ -\xi(s)f_N'(s)+\ff {g_1(s)} 2 f_N''(s) =-1,\ \ s\in (0,N].$$ Thus, by (\ref{B1}), $\d\<M\>_t\ge g_1(\|X_t\|_H^2)\d t$ and  It\^o's formula,
$$\d f_N(\|X_t\|_H^2)\le-\d t +f_N'(\|X_t\|_H^2)\d M_t,\ \ \ t<\tau_0\land\tau_N.$$ Therefore,
\beq\label{B3'} \E(\tau_0\land \tau_N) \le f_N(\|x\|_H^2) <\infty.\end{equation}
On the other hand, taking
$$f(s):=\int_0^s \exp\bigg[2\int_1^t\ff{\xi(u)}{g(u)}\d u\bigg]\d t,\ \ s\ge 0,$$ we have
$$-\xi(s)f'(s) +\ff{g(s)} 2 f''(s) =0,\ \ \ s>0.$$ Moreover, by assumption and the definition of $g$
 $$f''(\|X_t\|_H^2)\d\<M\>_t \le f''(\|X_t\|_H^2)g(\|X_t\|_H^2)\d t.$$ Therefore, by (\ref{B1}), $f(\|X_t\|_H^2)$ is a super-martingale up to time $\tau_0\land\tau_N$. So,
$$f(\|x\|_H^2)\ge \E f(\|X_{t\land\tau_0\land\tau_N}\|_H^2)\ge \P(\tau_N\le t\land\tau_0)f(N),\ \ \ t>0.$$ Combining this with (\ref{B3'}) we obtain
$$\P(\tau_0=\infty) = \P(\tau_0=\infty, \tau_N\land\tau_0<\infty)\le\lim_{t\to\infty} \P(\tau_N\le \tau_0\land t)\le \ff{f(\|x\|_H^2)}{f(N)}.$$ Then the first part of the assertion follows  by letting $N\to\infty$. The second follows by realizing that the denominator of the right-hand side of (\ref{C'}) is equal to infinity if $g(s)\ge \rr_3 s^2$ for large $s$, and that for large $s$ one has $\xi(s)<0$ and hence $g(s)=g_1(s).$
\end{proof}

\section{Applications to stochastic fast-diffusion equations}

The aim of this section is to apply Theorems \ref{T2.1} and \ref{T2.2} to a class of SPDEs as mentioned in Remark 2.2
and in the Introduction.

Let $(E,\scr E,\mu), L$ and $P_t$ be as  in Section 2 such that the spectrum of $L$ is  discrete with strictly positive eigenvalues $\{\ll_k\}_{k\ge 1}$ counting multiplicities with increasing order, and let $\{e_k\}_{k\ge 1}$ be the corresponding unit eigenfunctions forming an ONB of $L^2(\mu)$. Then $H$ is the completion of $L^2(\mu)$ w.r.t. the inner product
$$\<x,y\>_H:=\<x, L^{-1} y\>_2 =\sum_{k=1}^\infty \ff{\<x, e_k\>_2\<y, e_k\>_2}{\ll_k},$$ where $\<\cdot,\cdot\>_2$ also denotes both the inner product in   $L^2(\mu)$ and its extension to $H\times H^*$. Next, let
\beq\label{U} \|P_t\|_{1\to \infty}\le (c_\infty t)^{-d/2},\ \ t>0\end{equation} hold for some constants $d,c_\infty\in (0,\infty).$ By (\ref{U}), $\|P_t\|_{2\to 2}\le \e^{-\ll_1t}$, and the Riesz-Thorin interpolation theorem, we have
$$\|P_t\|_{1+r\to \ff{1+r}r}\le (c_\infty t)^{-\ff{d(1-r)}{2(1+r)}}\e^{-\ll_1 (1-r)t/(1+r)},\ \ \ t>0.$$ Therefore,
\beq\label{GG} \gg(\theta)<\infty,\ \ \text{if}\ \  \theta\in \Big(0, \ff{2(1+r)}{d(1-r)}\Big)\cap (0,1].\end{equation}
A standard example for the framework is that $L=-\DD+V$ for the Dirichlet Laplacian $\DD$ on a   domain $\O\subset \R^d$ having finite volume, and for a nonnegative locally bounded measurable function $V$ on $\O$. It is also the case if $\O$ has infinite volume but
$$\mu(V\le r):= \mu(\{x\in \O: V(x)\le r\})<\infty,\ \ r\ge 0,$$ where $\mu$ stands for the Lebesgue measure on $\O$. In this case, according to \cite{WW}, $L$  has discrete spectrum as well.

Moreover, let $\Psi: \R\to 2^\R$ be a maximal monotone graph such that $0\in \Psi(0)$ and
\beq\label{P} C(  |s|^{q}+ |s|)\ge sy\ge \theta_1 |s|^{1+r}+\theta_2 s^{2},\ \ s\in\R, y\in\Psi(s)\end{equation} holds for
some constants $r\in [0,1), C>0, \theta_1>0, \theta_2\ge 0, q\ge 1+r.$

Now, let $W_t$ be a cylindrical Brownian motion on $L^2(\mu)$.
Extending   the framework investigated in \cite{BDRa, BDRb, BDPR11}, where $L=-\DD$ on a bounded domain in $\R^d$ (
for $d\le 3$),
 we  consider the following SPDE on $H$:
\beq\label{EE} \d X_t + L\Psi(X_t)\d t \ni   B(X_t)\d W_t,\end{equation} where $B: H\to \scr L_2(L^2(\mu);H)$ is measurable,
subject to conditions to be specified in the following   Subsections 3.1-3.3.
\

Following Definition 2.1 in \cite{BDRb}  we call  a continuous adapted process $X:=(X_t)_{t\ge 0}$ on $H$ a solution to (\ref{EE}) if:    \beg{enumerate}\item[(a)] $(X_t)_{t\in [0,T]}\in L^{2}([0,T]\times \OO\to L^{2}(\mu); \d t\times \P)$ for any $T>0$.
%\cap L^2([0,T]\times\OO\to H; \d t\times \P).$ If $\rr_2>0$ then $(X_t)_{t\in [0,T]}\in L^{1+p}([0,T]\times \OO\to L^{1+p}(\mu); \d t\times \P)$ also holds;
\item[(b)] There exists a progressively measurable  process $\eta:=(\eta_t)_{t\ge 0}$ such that $ (\eta_t)_{t\in [0,T]}\in L^2([0,T]\times\OO\to L^2(\mu);\d t\times\P)$ for any $T>0$, $\eta\in \Psi(X)\ \d t\times\P\times\mu$-a.e., and   $\P$-a.s.
$$\<X_t, e_k\>_H=\<X_0,e_k\>_H - \int_0^t \<\eta_s,  e_k\>_2\d s +\int_0^t \<B(X_s)\d W_s, e_k\>_H,\ \ t\ge 0, k\ge 1.$$\end{enumerate}
Applying It\^o's formula to $\<X_t,e_k\>_H^2,$ multiplying by $\ll_k$ and summing over $k\in\mathbb N$, it follows that
  \beq\label{IT}\d\|X_t\|_H^2 = -2\<X_t, \eta_t\>_2\d t +\|B(X_t)\|_{\scr L_2(L^2(\mu);H)}^2\d t +2\<B(X_t)\d W_t, X_t\>_H.\end{equation}

\subsection{Extinction   for stochastic fast-diffusion type equations with linear multiplicative noise}

In this subsection we consider  $L=(-\DD+V)^\aa$ for $\DD$ the Dirichlet Laplacian on a domain $\O\subset \R^n$ and $V\ge 0$ being a measurable function on $\O$ such that the spectrum of $L$ is discrete, where $n\in \mathbb N$ and $\aa\in (0,1]$  are fixed constants. Let $\mu$ be the Lebesgue measure on $\O$. It is well-known that the   semigroup $P_t^{(0)}:=\e^{t(\DD-V)}$ satisfies
$$\|P_t^{(0)}\|_{1\to\infty} \le c_0t^{-n/2},\ \ t>0$$ for some constant $c_0\in (0,\infty)$. Then
(\ref{U}) holds for $P_t:=\e^{-tL}$ with $d=\ff n \aa$ and some constant $c_\infty\in (0,\infty).$
Consider (\ref{EE}) for
\beq\label{BB1'}B(x)h=\sum_{k=1}^\infty \mu_k \<h, e_k\>_2x e_k,\ \ x\in L^{2}(\mu), h\in L^2(\mu),\end{equation} where $\{\mu_k\}_{k\ge 1}$ is a sequence of constants such that
\beq\label{T} \sum_{k=1}^\infty \mu_k^2 \ll_k^{\ff d 2\lor (1+\vv)}<\infty\end{equation} holds for some constant $\vv>0.$
When $L=-\DD$ for $\DD$ the Dirichlet Laplacian on a  bounded  domain in $\R^d$,
$$\ll_k =\text{O}( k^{2/d}) $$ holds for large $k$,  so that (\ref{T}) follows from
$$\sum_{k=1}^\infty \mu_k^2k^{1\lor(\ff 2 d +\vv)}<\infty$$ for some $\vv>0.$  We note that when $d\ne 2$ (\ref{T}) is weaker than the corresponding condition  used in \cite{BDPR11} (in \cite{BDRa, BDRb} the condition is even stronger as only finite modes of the noise are allowed).
%, where the  power of $\ll_k$ is   $1\lor (d-1)$ rather than $\ff d 2$ (but when $d=2$ they are equal).

To ensure that   $(H1)$ holds for  solutions to  (\ref{EE}), we first observe that (\ref{T}) implies $\|B(x)\|_{\scr L_2(L^2(\mu);H)}^2\le \rr_3\|x\|_H^2$ for some constant $\rr_3>0.$

\beg{prp}\label{P3.1}  Let $B$ be as in $(\ref{BB1'})$. If $(\ref{T})$ holds, then the linear map
 $L^2(\mu)\ni x\mapsto B(x)\in \scr L_2(L^2(\mu), H)$ extends by  continuity to all of $H$ and
\beq\label{R} \rr_3:= \ff 1 2 \sup_{\|x\|_H^2=1} \|B(x)\|_{\scr L_2(L^2(\mu);H)}^2<\infty,\end{equation} so that $\|B(x)\|_{\scr L_2(L^2(\mu);H)}^2\le 2\rr_3\|x\|_H^2$ holds for all $x\in H$.\end{prp}

To prove this result, we need the following lemma.

\beg{lem}\label{LL} Let $(\EE,\D(\EE))$ be a symmetric Dirichelt form on $L^2(\mu)$ over a $\si$-finite measure space $(E,\F,\mu)$,
and let $(\scr L,\D(\scr L))$ be the associated Dirichlet operator, i.e. the negative definite self-adjoint operator on $L^2(\mu)$ associated to the symmetric form $(\EE,\D(\EE))$.
Then for any $g\in L^\infty(\mu)\cap\D(\scr L)$ such that $Lg\in L^\infty(\mu)$,
$$\EE(fg, fg)\le \|g\|_\infty^2\EE(f,f)- \mu(f^2g \scr L g),\ \ \ f\in\D(\EE).$$\end{lem}
\beg{proof} For $t>0$, let $J_t$ be the symmetric measure on $E\times E$ such that
$$J_t(A\times B)=\mu(1_A P_t 1_B),\ \ \ A,B\in \F,$$ where $P_t$ is the associated Markov semigroup. Then, by the symmetry of $J_t$,
\beg{equation*}\beg{split} \EE(h,h)&=\lim_{t\downarrow 0} \ff 1 t \int_E h(h-P_t h)\d\mu = \lim_{t\downarrow 0} \ff 1 t  \int_{E\times E} h(x) \Big(\ff {h} {P_t 1}(x)-h(y)\Big)J_t(\d x,\d y) \\
&=\lim_{t\downarrow 0} \ff 1 t \bigg(\int_{E\times E} h(x)\big(h(x)-h(y)\big)J_t(\d x,\d y) +
\int_T h^2 (  1 -P_t 1)\d\mu\bigg)\\
&=\lim_{t\downarrow 0} \ff 1 t \bigg(\ff 1 2 \int_{E\times E} \Big\{h(x)-h(y)\Big\}^2J_t(\d x,\d y)+
\int_T h^2  ( 1 - P_t 1)\d\mu\bigg),\ \ h\in\D(\EE).\end{split}\end{equation*}  Combining this with
$$\Big\{(fg)(x)-(fg)(y)\Big\}^2 = (f(x)-f(y))^2 g(x)g(y) +\big(f(y)^2g(y) -f(x)^2 g(x)\big)\big(g(y)-g(x)\big),  $$ and using the symmetry of $J_t$, we arrive at
\beg{equation*}\beg{split} &\EE(fg,fg) = \lim_{t\downarrow 0} \ff 1 t \bigg(\ff 1 2 \int_{E\times E} \Big\{(fg)(x)-(fg)(y)\Big\}^2J_t(\d x,\d y)+
\int_T f^2g^2  ( 1 -P_t 1)\d\mu\bigg)\\
&\le  \lim_{t\downarrow 0} \ff 1 t \int_{E\times E} \Big(\ff {\|g\|_\infty^2} {2} \big\{f(x)-f(y)\big\}^2   -   f(x)^2g(x) \big(g(y)-g(x)\big)\Big)J_t(\d x,\d y)\\
&= \lim_{t\downarrow 0}\ff 1 t \bigg\{\ff{\|g\|_\infty^2}{2} \int_{E\times E}   \Big\{f(x)-f(y)\Big\}^2 J_t(\d x,\d y)
 - \int_E \Big\{f^2 g(P_t g-g)+ f^2g^2 (1-P_t 1)\Big\}\d\mu\bigg\}\\
&\le  \|g\|_\infty^2 \EE(f,f)- \mu(f^2g\scr L g).\end{split}\end{equation*}
\end{proof}

\beg{proof}[Proof of Proposition \ref{P3.1}]  First we note that
by (\ref{U}),
\beq\label{R3.8}\|e_k\|_\infty =\e \|P_{1/\ll_k}e_k\|_\infty\le \e \|P_{1/\ll_k}\|_{2\to\infty}\le c'   \ll_k^{d/4}\end{equation} holds for some constant $c'>0.$
Since $B(x)$ is linear in $x$, it suffices to prove (\ref{R}). By the definition of $B(x)$, we have
\beq\label{B} \|B(x)\|_{\scr L_2(L^2(\mu);H)}^2=\sum_{k=1}^\infty \|B(x)e_k\|_H^2 =\sum_{k=1}^\infty \mu_k^2\|xe_k\|_H^2.\end{equation} So, it suffices to show that   there exists a constant $c>0$ such that
\beq\label{B2} \|x e_k\|_H^2 \le c \|x\|_H^2 \ll_k^{\ff d 2\lor (1+\vv)},\ \ x\in H, k\ge 1.\end{equation} By (\ref{R3.8}) and an approximation argument, we only need to prove this inequality for $x\in L^2(\mu)$. In this case
\beq\label{BB1}\|xe_k\|_H=\sup_{\EE(f,f)\le 1} \mu(xe_k f) \le \|x\|_H\sup_{\EE(f,f)\le 1} \ss{\EE(e_k f, e_k f)},\end{equation}
where $(\EE,\D(\EE))$ is the associated Dirichlet form associated to  $\scr L:=-L$.
%Next, by (\ref{U}), the super Poincar\'e inequality
%\beq\label{SP} \mu(g^2)\le r \EE(g,g) + c_0(1+r^{-d/2})\mu(|g|)^2,\ \ \ r>0, g\in \D(\EE)\end{equation} holds for some constant $c_0>0$, see
%\cite[Theorem 4.5 (2)]{W00} or \cite[Theorem 3.3.15(2)]{Wbook}.  Then
%$$\mu(f^2 e_k^2)\le \ff 1 {2\ll_k} \EE(e_kf, e_kf) +c_0\big(1+(4\ll_k)^{d/2}\big)\mu(|e_kf|)^2\le \ff 1 {4\ll_k}
%\EE(e_kf, e_kf) +c_1\ll_k^{d/2}\mu(f^2)$$ holds for some constant $c_1>0$ as $\mu(e_k^2)=1$ and $\ll_k\ge \ll_1>0$.
 Since $\<f^2e_k, L e_k\>_2=\ll_k \mu(f^2e_k^2)$, by Lemma \ref{LL} with $\scr L=-L$,
we obtain \beq\label{WK} \EE(e_k f, e_k f)\le  \|e_k\|_\infty^2 +\ll_k \mu(e_k^2f^2).\end{equation}
For $p:=d\lor (2+\ff \vv d)$, it follows from (\ref{U}) and $\|P_t\|_2\le \e^{-\ll_1 t}$ that
$$\|P_t\|_{1\to\infty}\le c t^{-p/2},\ \ t>0$$ holds for some constant $c>0.$ This implies the Sobolev inequality (cf. \cite{D})
$$\|f\|_{\ff{2p} {p-2}}^2\le C\EE(f,f),\ \ f\in\D(\EE) $$ for some constant $C>0$.
Combining this  with (\ref{WK}) and noting that $\|e_k\|_2=1$, we obtain
\beq\label{WFY}\EE(e_kf,e_kf)\le \|e_k\|_\infty^2 + C\ll_k \|e_k\|_p^2\le \|e_k\|_\infty^2+C\ll_k\|e_k\|_\infty^{\ff{2(p-2)}p}.\end{equation} Combining (\ref{R3.8})  with (\ref{BB1}) and (\ref{WFY}) we prove (\ref{B2}). \end{proof}

Now exactly the same arguments as in the proof of \cite[Theorem 2.2]{BDRb} imply that for any
$x\in L^{4\lor (2q)}(\mu)$, there exists a unique solution to (\ref{EE}) with $X_0=x$, which is non-negative if $x\ge 0$. Hence
by Proposition \ref{P3.1} and (\ref{IT}), we have the following consequence of Theorem \ref{T2.1}:

\beg{cor}\label{C1} Assume $X_0=x\in L^{4\lor 2q}(\O)$ and that $(\ref{P})$  and $(\ref{T})$ hold. Let $\rr_1=\theta_1,\rr_2=\theta_2$,
 and let $\rr_3$ be as in Proposition $\ref{P3.1}$. Then for any
$\theta\in \big(0, \ff{2(1+r)}{d(1-r)}\big)\cap (0,1], d:= \ff n \aa$, all assertions in Theorem $\ref{T2.1}$ hold for
solutions to $(\ref{EE})$ with $B$ given in $(\ref{BB1}).$ \end{cor}

\paragraph{Example 3.1.}
{\bf (i) (SOC-case)}. Let $\Psi$ be as in (\ref{R1.1'}). In this case (\ref{P}) holds with $r=0, $
some $\theta_1>0, C>0$ and $q=1, \theta_2=0$ in the BTW-model and $q=2,\theta_2>0$ in the Zhang-model respectively for
$s\ge 0.$ But if $X_0=x\in L^4(\O), x\ge 0$, as mentioned above, we have $X_t\ge 0$ for all $t\ge 0.$ Therefore,
 to consider $s\ge 0$ is sufficient, since we may change $\Psi$ on $(-\infty,0)$ to become an odd function without changing anything.
 Hence we can apply Proposition \ref{P3.1} with $B$ as in (\ref{BB1'}) satisfying (\ref{T}). In the BTW-model $\theta_2=0$,
 so we need $\ff{2\aa}n>1$, i.e. $n=1$ and $\aa>\ff 1 2,$ and we have only extinction in finite time with positive probability
 for small enough noise or initial conditions $X_0=x\in L^4(\O), x\ge 0,$ with small enough $H$-norm.  So, we recover the corresponding result from \cite{BDRb, BDPR11} in the special case $\aa=1$. Furthermore, if $B=0$, hence $\rr_3=0$ in Theorem \ref{T2.1}(3), which thus applies for $n=1,\aa>\ff 1 2,$ to give extinction, recovering the deterministic case from \cite{BDRb, BDPR11}. In the Zhang-model, however,
 we have $\rr_2=\theta_2>0$. Hence for all $\aa\in (0,1]$ and all dimensions $n$, if we choose $\theta\in (0,\ff {2\aa}n),$  and
 for small enough noise we can apply Theorem \ref{T2.1}(2) to get exponential integrability of the extinction time, hence in particular extinction in finite time with probability one, provided $X_0=x\in L^4(\O), x\ge 0.$  In particular, for the deterministic Zhang-model we have $r=\rr_3=0$ and $\rr_1=\theta_1,\rr_2=\theta_2>0$, so that Theorem \ref{T2.1}(2)   implies that
$$ \tau_0=\lim_{\bb\to 0} \ff 1 \bb (\e^{\bb \tau_0}-1)\le \inf_{\theta\in (0,1]} \ff{\gg(\theta)^{\ff 1 2 } \|X_0\|_H^\theta}{\theta \theta_1^\theta\theta_2^{1-\theta}\ll_1^{\ff{1-\theta}2}}<\infty.$$

{\bf  (ii) (Fast diffusion-case)}. Let $\Psi$ be as in (\ref{R1.4'}). Then (\ref{P}) holds with $\theta_1=1, C=1, q=1+r, \theta_2=0.$ So,
 if $B$ is not identically equal to zero, only Theorem \ref{T2.1}(1) applies with $\rr_2=\theta_2=0$ and $\theta=1$. So,
 we need $\ff{2\aa(1+r)}{n(1-r)}>1$, i.e. $r>\ff{n-2\aa}{n+2\aa}$, and we only have extinction in
  finite time with positive probability for small enough noise or initial conditions with small enough $H$-norm.
  So, again we recover the corresponding results from \cite{BDRa, BDPR11} in the special case $\aa=1, n\le 3$ (though for $n=3,\aa=1,$ also the case $r=\ff{n-2\aa}{n+2\aa}+\ff 1 5$ is covered in \cite{BDRa, BDPR11}). If $B=0$,
  Theorem \ref{T2.1}(3) applies with $\theta=1$, leading to the same restriction $r>\ff{n-2\aa}{n+2\aa}.$
  Hence we get extinction in finite time for the deterministic case, which appears, however,
  to be a new result if $\aa<1$. For $\aa =1$   we recover the well-known results for the
  deterministic  fast diffusion equation. Finally, we point out, that, adding a linear term to $\Psi$ we again get extinction in finite time with probability one in the same way as in the Zhang-model above.

 \

In the next subsection, we consider   much stronger noises such that $\mu_k=1$ in (\ref{BB1'}) is allowed.

\subsection{Exponential integrability of $\tau_0$   for strongly dissipative equations with uncoloured linear multiplicative noise}

 In addition to (\ref{P}), we assume that
\beq\label{P'} \Psi\in C(\R),\ \
(s-t)(\Psi(s)-\Psi(t))\ge \kk |s-t|^{2},\ \ \ \ \ s,t\in\R \end{equation} holds for some constants $\kk>0.$ Because $\Psi(0)=0$, (\ref{P'}) implies that
(\ref{P}) holds with $\theta_2>0$.
 A simple example of $\Psi$ satisfying (\ref{P}) and (\ref{P'}) is
$$\Psi(s)= \theta_1 |s|^r\text{sgn}(s)+\theta_2s$$ for some constants $r\in(0,1), \theta_1,\theta_2>0.$
Let $B_0$ be a bounded linear operator on $L^2(\mu)$.  Let
$$\rr_0:=  \ff 1 2 \|B_0\|_{2\to 2}^2 \sum_{k=1}^\infty \ff {\|e_k\|_\infty^2} {\ll_k}\in[0,\infty],$$ where $\|B_0\|_{2\to 2}$ is
the operator norm of $B_0$ in $L^2(\mu)$.
We consider the following stochastic differential equation on $H$:
\beq\label{E3} \d X_t = -L \Psi(X_t)\d t+ X_t B_0\d W_t.\end{equation}
It will turn out that if $\rr_0<\infty$, then this equation is a special case of  (\ref{EE}), since in this case the operator-valued function $B$ defined by
$$B(x)h= x B_0 h,\ \ \ h\in L^2(\mu),\ x\in L^\infty(\mu)$$ extends by continuity to all $x\in L^2(\mu)$ (see (\ref{3.15'}) below) which is sufficient under condition (\ref{P'}).
We emphasize here that $\rr_0<\infty$, e.g. if $L$ is the Dirichlet Laplacian on $(0,1)$ (cf. Remark 3.1 below).
Since $\Psi$ is single-valued we can use the result on uniqueness and existence of solutions from \cite{RRW}, which holds even for random initial conditions, as then does   the following theorem.

\beg{thm} \label{TTT} Assume that $(\ref{U}), (\ref{P}), (\ref{P'})$  hold. If  $\rr_0\in (0, \kk]\cap (0, \theta_2)$, then for any $X_0 \in L^2(\OO\to H,\F_0;\P)$, the equation $(\ref{E3})$ has a unique solution in the sense of \cite{RRW}. If moreover $(\ref{U})$ holds, then for any  $\theta\in \big(0, \ff{2(1+r)}{d(1-r)}\big)\cap (0,1]$ and any $\bb\in \big(0, \theta (1-r)(\theta_2-\rr_0)\ll_1\big),$
$$\E \e^{\bb\tau_0}\le 1+\ff{\bb\gg(\theta)^{(1+r)/2}\E\|X_0\|_H^{\theta(1-r)}}{\theta (1-r)\theta_1^\theta\{\theta_2-\bb/(\ll_1\theta(1-r)\}^{1-\theta}\ll_1^{(1-r)(1-\theta)/2}}<\infty$$ provided $\E\|X_0\|_H^{\theta(1-r)}<\infty.$
\end{thm}

\beg{proof} For the existence and uniqueness of solutions, we need only to verify the assumptions of \cite[Theorem 2.1]{RRW}. To this end, we take ${\bf V}=L^2(\mu)$, $A(x)= - L\Psi(x), R(x)=\mu(x^2)$ and for fixed $T>0$, $K= L^2([0,T]\times\OO\times E; \d t\times \P\times\mu).$ Then  assumptions {\bf (K)},  {\bf (H1)} and {\bf (H4)} in \cite{RRW} follow immediately from (\ref{P}) and the continuity of $\Psi$. It remains to verify that for some constants $c_1,c_2,c_4\in\R$ and $c_3>0$
\beg{enumerate}\item[{\bf (H2)}] $-2\<\Psi(x)-\Psi(y), x-y\>_2+\|B(x)-B(y)\|_{\scr L_2(L^2(\mu);H)}^2\le c_1 \|x-y\|_H^2,\ \ x,y\in L^2(\mu);$
\item[{\bf (H3)}]   $-2\<\Psi(x),x\>_2+ \|B(x)\|_{\scr L_2(L^2(\mu);H)}^2\le c_2\|x\|_H^2-c_3\|x\|_{2}^{2}+c_4,\ \ x\in L^2(\mu).$\end{enumerate}
Observe that for $x\in L^\infty(\mu)$
\beg{equation}\label{3.15'}\beg{split} \|B(x)\|_{\scr L_2(L^2(\mu);H)}^2&= \sum_{k=1}^\infty\|xB_0 e_k\|_H^2 =\sum_{k,j=1}^\infty \ff{\mu(\{B_0e_k\}e_jx)^2}{\ll_j}=\sum_{k,j=1}^\infty \ff{\mu(e_k B_0^*(e_jx))^2}{\ll_j}\\
&=\sum_{j=1}^\infty \ff{\|B_0^*(e_jx)\|_2^2}{\ll_j}
 \le  \|B_0\|_{2\to 2}^2\|x\|_{2}^2 \sum_{j=1}^\infty \ff{\|e_j\|_\infty^2}{\ll_j}= 2\rr_0\|x\|_{2}^2.\end{split}\end{equation}
 Since $B(x)-B(y)= B(x-y)$,
combining (\ref{3.15'}) with (\ref{P'}) and noting that $\kk\ge \rr_0$, we obtain  {\bf (H2)} for $c_1=0$. Moreover,  (\ref{P}) and (\ref{3.15'})
imply that {\bf (H3)} holds for $c_2=c_4=0$ and $ c_3= \theta_2-\rr_0>0$.
Next, by  (\ref{IT}), we have
$$\d \|X_t\|_H^2\le -2\big\{\theta_1\|X_t\|_{1+r}^{1+r}+ \theta_2 \|X_t\|_{2}^{2}-\rr_0\|X_t\|_{2}^2\big\}\d t + \d M_t,$$ where $$\d M_t = 2 \< X_t B_0\d W_t, X_t\>_H.$$
This implies $(H1)$ for $\rr_1=\theta_1,\rr_2=\theta_2-\rr_0>0$ and $\rr_3=0.$  Then the desired result on $\tau_0$ follows from Theorem \ref{T2.1}(2).
\end{proof}
\paragraph{Remark 3.1.} Let e.g. $L=-\DD$ on $(0,1).$ Then  we have
$$ \ll_k= \pi^2k^2, \ \ \ e_k(s)= \ss 2 \sin(\pi k s),\ \ k\ge 1, s\in (0,1),$$ so that $\rr_0<\infty.$
Therefore, $\rr_0<\infty$ and Theorem \ref{TTT} applies, for small enough $\|B_0\|_{2\to 2}$.

\

In the next subsection, we consider the case with a   noise   having a component   in the direction $X_t$, so that Theorem \ref{T2.2} applies.

\subsection{Extinction with probability $1$ for special noise}

Now let us consider again the situation of Subsection 3.1. Let $e$ be an unit element in $L^2(\mu)$ and let
$$B(x)h= cx\<h,e\>_2,\ \ x\in H, h\in L^2(\mu),$$ where $c\ne 0$ is a constant.   Taking $B_t:=\<W_t,e\>_2$ (a one-dimensional Brownian motion)   (\ref{EE}) reduces to
\beq\label{EE2} \d X_t +L\Psi(X_t)\ni cX_t \d B_t.\end{equation}  By (\ref{P}) and  (\ref{IT}), $(H1)$ holds for $\rr_1=\theta_1, \rr_2=\theta_2$ and
$$\d M_t= 2c\<X_t, X_t\>_H \d t.$$ We have
$$\d\<M\>_t= 4 c^2\|X_t\|_H^4\d t=g(\|X_t\|_H^2)\d t$$ for $g(s)= 4 c^2s^2.$

\beg{cor}\label{C2} Assume that $(\ref{U})$ and  $(\ref{P})$ hold, and let $\rr_3=\ff {c^2} 2.$   If either $\theta_2>0$ or $\theta_2=0$ but $\ff{2(1+r)}{d(1-r)}>1$, then $\P(\tau_0<\infty)=1$ holds for
solutions to $(\ref{EE2})$ for any $X_0=x\in H$. \end{cor}
\beg{proof} Obviously, $g(s)\ge 2\rr_3 s^2$ and $(H1)$ hold. By Theorem \ref{T2.2} for $g_1=g_2=g$, it suffices to note that due to (\ref{GG}) one has
$\gg(\theta)<\infty$ for all $\theta\in \big(0,\ff{2(1+r)}{d(1-r)}\big)\cap (0,1]$, and thus $\gg(1)<\infty$ if $\ff{2(1+r)}{d(1-r)}>1.$\end{proof}

\paragraph{Remark 3.2.} We note that to make stochastic perturbations in directions other than $X_t$, we may consider
 \beq\label{EE2} \d X_t +L\Psi(X_t)\ni c X_t \d B_t+\tt B(X_t)\d\tt W_t,\end{equation}where  $\tt W_t$ be a cylindrical Brownian motion on $L^2(\mu)$ which is independent of $W_t$, and $\tt B: H\to \scr L_2(L^2(\mu);H)$ is such that $\|B(x)\|_{\scr L_2(L^2(\mu);H)}^2\le \tt c^2\|x\|_H^2.$
 Then Theorem \ref{T2.1} applies to $\rr_3=\ff 1 2 (c^2+\tt c^2)$, while the assertion in Theorem \ref{T2.2} holds for this $\rr_3$ and $g_1(s)= 4c^2s^2, g_2(s)= 4(c^2+\tt c^2)s^2.$

 \

Finally, we consider one more case which in fact generalizes the one considered above (take $N=0$ below), but  for, which the noise exists not only in one, but  in finitely many directions, and both Theorem \ref{T2.1} and Theorem \ref{T2.2} apply:  Let $N\in\N$ and consider (\ref{EE}) for
\beq\label{BB}B(x)h= \sum_{k=1}^N \mu_k \<x, e_k\>_2\<h,e_k\>_2e_k + \mu_{N+1} \<h,e_{N+1}\>_2 \pi_N^\bot(x),\ \ \ x,h\in H,\end{equation}
 where $\{\mu_k\}_{k= 1}^{N+1}\subset\R$ and
 $$\pi_N^\bot(x):= \sum_{k=N+1}^\infty \<x,e_k\>_2 e_k.$$  Let
 $$\rr_3= \ff 1 2 \sup_{1\le k\le N+1}\mu_k^2.$$ We have
$$\|B(x)\|_{\scr L_2(L^2(\mu);H)}^2 =\sum_{i=1}^{N +1}\|B(x)e_k\|_H^2=\sum_{k=1}^{N} \mu_k^2 \ff{\<x,e_k\>_2^2}{\ll_k}+\mu_{N+1}^2\|\pi_N^\bot(X_t)\|_H^2\le 2\rr_3\|x\|_H^2.$$
By  (\ref{P}) and (\ref{IT}),  $(H1)$ holds for $\rr_1=\theta_1, \rr_2=\theta_2$ and
$$\d M_t = 2\<B(X_t)\d W_t, X_t\>_H= 2\sum_{k=1}^N \ff{\mu_k}{\ll_k}\<X_t,e_k\>_2^2\<\d W_t, e_k\>_2+ 2\mu_{N+1}\|\pi_N^\bot(X_t)\|_H^2\<\d W_t, e_{N+1}\>.$$ Obviously,
\beg{equation*}\beg{split} &4\Big(\sum_{k=1}^{N+1} \mu_k^2\Big)\|X_t\|_H^4 \d t\ge \d\<M\>_t=4\Big\{\sum_{k=1}^N \ff{\mu_k^2}{\ll_k^2} \<X_t,e_k\>_2^4+\mu_{N+1}^2 \|\pi_N^\bot(X_t)\|_H^4\Big\}\d t\\
&\ge \ff{4\big(\sum_{k=1}^N \ff{\<X_t,e_k\>_2^2}{\ll_k}+\|\pi_N^\bot(X_t)\|_H^2\big)^2}{\sum_{k=1}^{N+1} \mu_k^{-2}}\,\d t=\ff{4\|X_t\|_H^4\d t}{\sum_{k=1}^{N+1}\mu_k^{-2}}.
 \end{split}\end{equation*} Let
$$c_1=\ff 4 {\sum_{k=1}^{N+1}\mu_k^{-2}},\ \ c_2= 4\sum_{k=1}^{N+1} \mu_k^2.$$ Therefore, if $\mu_k^2>0$ for $1\le k\le N+1$, then (\ref{C}) holds for $g_i(s)= c_is^2,\ i=1,2.$

\beg{cor}  Assume that $(\ref{U})$ and $ (\ref{P})$  hold. Let $\rr_3= \ff 1 2 \sup_{1\le k\le N+1}\mu_k^2.$  Then for any
$\theta\in \big(0, \ff{2(1+r)}{d(1-r)}\big)\cap (0,1]$, all assertions in Theorems $\ref{T2.1}$   hold for
solutions to $(\ref{EE})$ with $B$ given in $(\ref{BB}).$  If moreover $\mu_k^2>0$ for $1\le k\le N+1$, the assertion in Theorem $\ref{T2.2}$ holds for $g_i(s)=c_is^2,\ i=1,2.$\end{cor}

\beg{thebibliography}{99}

\bibitem{BTW88} P. Bak, C. Tang, K. Wiesenfeld, \emph{Self-organized criticality,} Physical Review A, 38(1988), 3611--374.

\bibitem{BJ92} P. B\'antay, I. M. J\'anosi, \emph{Self-organized criticality and anomalous diffusion,} Physica A 185(1992), 11--18.

\bibitem{BDR08} V. Barbu, G. Da Prato, M. R\"ockner, \emph{Existence and uniqueness of nonnegative solution to the stochastic porous media equation,} Indiana Univ. Math. J. 57(2008), 187--211.
\bibitem{BDRa}  V. Barbu, G. Da Prato, M. R\"ockner, \emph{Finite time extinction for solution to fast diffusion stochastic porous media equations,} C. R. Acad. Sci. Paris-Math. 347(2009), 81--84.
\bibitem{BDRb}  V. Barbu, G. Da Prato, M. R\"ockner, \emph{Stochastic porous media equation and self-organized criticality,} Comm. Math. Phys. 285(2009), 901--923.

\bibitem{BDPR11}  V. Barbu, G. Da Prato, M. R\"ockner, \emph{Finite time extinction of solutions to fast diffusion equations driven by linear multiplicative noise,} CRC701-Preprint, 29pp. (2011).

\bibitem{D}  E. B. Davies, \emph{Heat Kernels and Spectral Theory,}
Cambridge: Cambridge Univ. Press, 1989.

\bibitem{MR92} Z.-M. Ma  and M. R\"ockner, \emph{Introduction to
the Theory of (Non-Symmetric) Dirichlet Forms,} Springer-Verlag,
Berlin, 1992.

\bibitem{RRW} J. Ren, M. R\"ockner, F.-Y. Wang, \emph{Stochastic generalized porous media and fast-diffusion equations, } J. Diff. Equations,  238(2007), 118--152.

\bibitem{RW08}	 M. R\"ockner, F.-Y. Wang, \emph{ Non-monotone stochastic generalized porous media equations,} J. Differential Equations 245(2008), 3898-3935.

\bibitem{WW}	F.-Y. Wang, J.-L. Wu, \emph{Compactness of Schr\"odinger semigroups with unbounded below potentials,} Bulletin des Sciences Mathematiques, 132(2008), 679--689.

\bibitem{V10} A. de Pablo, F. Quir\'os, A. Rodriguez, J. L. V\'azquez, \emph{A fractional porous medium equation,} Adv. Math. 226(2011), 1378--1409.

\bibitem{Z89} Y.C. Zhang,
\emph{Scaling theory of self-organized criticality,}
Phys. Rev. Lett. 63(1989), 470--473.
%\bibitem{W00} F.-Y. Wang,  \emph{Functional inequalities, semigroup properties and spectrum estimates,}  Infinite Dimensional Analysis, Quantum Probability and Related Topics 3:2(2000), 263--295.

%\bibitem{Wbook} F.-Y. Wang, \emph{Functional Inequalities, Markov Semigroups and Spectral Theory,} Science Press, Beijing, 2005.

 \end{thebibliography}

\end{document}